\documentclass[11pt,a4paper,reqno]{amsart}
\usepackage{amssymb,amsmath,amsthm}
\usepackage[utf8]{inputenc}
\usepackage[T1]{fontenc}
\usepackage{enumerate}
\usepackage[all]{xy}
\usepackage{fullpage}
\usepackage{comment}
\usepackage{array}
\usepackage{longtable}
\usepackage{stmaryrd}
\usepackage{mathrsfs} 
\usepackage{xcolor}
\usepackage{mathtools}

\def\gcd{\mathrm{gcd}}
\def\deg{\mathrm{deg}}

\usepackage{hyperref}
\hypersetup{
    colorlinks=true,
    linkcolor=blue,
    filecolor=red,  
    citecolor=green,
    urlcolor=cyan,
    pdftitle={GON},
    pdfpagemode=FullScreen,
    }

\usepackage{cleveref}
\crefformat{section}{\S#2#1#3}
\crefformat{subsection}{\S#2#1#3}
\crefformat{subsubsection}{\S#2#1#3}
\usepackage{enumitem}
\usepackage{tikz}
\usepackage{mathdots}

\newtheorem{theorem}{Theorem}[section]
\newtheorem{lemma}[theorem]{Lemma}

\newtheorem{corollary}[theorem]{Corollary}

\theoremstyle{definition}
\newtheorem{definition}[theorem]{Definition}

\theoremstyle{remark}
\newtheorem{remark}[theorem]{Remark}
\newtheorem{question}[theorem]{Question}

\makeatletter
\newcommand{\subalign}[1]{%
  \vcenter{%
    \Let@ \restore@math@cr \default@tag
    \baselineskip\fontdimen10 \scriptfont\tw@
    \advance\baselineskip\fontdimen12 \scriptfont\tw@
    \lineskip\thr@@\fontdimen8 \scriptfont\thr@@
    \lineskiplimit\lineskip
    \ialign{\hfil$\m@th\scriptstyle##$&$\m@th\scriptstyle{}##$\hfil\crcr
      #1\crcr
    }%
  }%
}
\makeatother

\numberwithin{equation}{section}

\title{Minimal Denominators Lying in Subsets of the Ring of Polynomials over a Finite Field}
\author{Noy Soffer Aranov}
\address{Graz University of Technology, Institute of Analysis and Number Theory, 8010 Graz, Austria}
\email{noy.sofferaranov@tugraz.at}
\begin{document}
\maketitle
\begin{abstract}
    Given a subset $\mathcal{S}\subseteq \mathbb{F}_q[x]$ and fixed integers $n,m\in \mathbb{N}$, we study the distribution of the smallest denominator $Q\in \mathcal{S}$ for which there exists $\mathbf{P}\in \mathbb{F}_q[x]^m$ such that $\left\Vert\frac{\mathbf{P}}{Q}-\boldsymbol{\alpha}\right\Vert<q^{-n}$, where $\boldsymbol{\alpha}\in x^{-1}\mathbb{F}_q((x^{-1}))^m$ is chosen randomly. We also consider the discrete analogue obtained by fixing a polynomial $N\in \mathbb{F}_q[x]$ with $\deg(N)=n$ and sampling $\boldsymbol{\alpha}$ uniformly from $\frac{1}{N}\mathbb{F}_q[x]^m$.
We prove that for any infinite subset $\mathcal{S}\subseteq \mathbb{F}_q[x]$, for every $n\in \mathbb{N}$ and every dimension $m$, the probability distributions of these two random variables coincide.
This result is significantly stronger than the corresponding statement in the real setting, where Balazard and Martin showed that the averages of the discrete and continuous smallest denominator functions are asymptotically close.
\end{abstract}
\section{Introduction}
Let $N\in \mathbb{N}$, and for $j\in \{0,1,\dots,N-1\}$, define the function 
$$q_j(N)=\min\left\{q\in \mathbb{N}:\exists p\in \mathbb{Z}:\frac{p}{q}\in \left[\frac{j}{N},\frac{j+1}{N}\right)\right\},$$
which returns the smallest denominator of a rational lying in an interval. Kruyswijk and Meijer \cite{KM} studied the sum $S(N)=\sum_{j=0}^{N-1}q_j(N)$ and proved that $S(N)\asymp N^{\frac{3}{2}}$. Moreover, they conjectured that $S(N)\sim \frac{16}{\pi^2}N^{\frac{3}{2}}$. Stewart \cite{St} improved the results of \cite{KM} by proving that $1.35N^{\frac{3}{2}}<S(N)<2.04N^{\frac{3}{2}}$. Balazard and Martin \cite{BM} proved the Kruyswijk Meijer conjecture and moreover showed that 
\begin{equation}
\label{eqn:BMAsymp}
S(N)=\frac{16}{\pi^2}N^{\frac{3}{2}}+O\left(N^{\frac{4}{3}}\ln^2N\right).\end{equation}
Shparlinski \cite{Sh} improved the error term in \eqref{eqn:BMAsymp} and showed that
$$S(N)=\frac{16}{\pi^2}N^{\frac{3}{2}}+O\left(N^{\frac{29}{22}+o(1)}\right).$$
To prove \eqref{eqn:BMAsymp}, Balazard and Martin compared the sum $S(N)$ and the average of the continuous analogue of the smallest denominator function. This analogue, which was proposed by Meiss and Sanders \cite{MS}, is defined as follows: For $\delta>0$ and $x\in [0,1)$, define 
$$q_{\min}(x,\delta)=\min\left\{q\in \mathbb{N}:\exists \frac{p}{q}\in \left(x-\frac{\delta}{2},x+\frac{\delta}{2}\right)\right\}.$$
Chen and Haynes \cite{CH} computed the distribution of the function $q_{\min}(x,\delta)$. The higher dimensional analogue of $q_{\min}(x,\delta)$ was studied by Artiles \cite{Art} and Marklof \cite{M} through dynamical methods. In particular, Chen and Haynes \cite{CH} (see also \cite{M} for a more general result in all dimensions) proved that
$$\int_{0}^1q_{\min}(x,\delta)dx=\frac{16}{\pi^2}\frac{1}{\delta^{\frac{1}{2}}}+O\left(\ln^2\delta\right).$$
To prove \eqref{eqn:BMAsymp}, Balazard and Martin \cite{BM} proved that
\begin{equation}
\label{eqn:R(N)}
S(N)-N\int_0^1q_{\min}\left(x,\frac{1}{N}\right)dx\ll N^{\frac{4}{3}}\ln^2 N.\end{equation}
In this paper, we extend the results of Balazard and Martin \cite{BM} to the function field setting in all dimensions and show that in this setting, the quantity analogous to the left hand side in equation \eqref{eqn:R(N)} is zero. Moreover, we prove that the continuous and discrete distributions are equal to each other, even when restricting the set of admissible denominators. First, the function field setting is introduced. 
\subsection{The Function Field Setting}
Let $q$ be a prime power, and denote the ring of polynomials over $\mathbb{F}_q$ by 
$$\mathcal{R}=\left\{\sum_{i=0}^na_ix^i:a_i\in \mathbb{F}_q,n\in \mathbb{N}\cup\{0\}\right\}.$$
For $n\in \mathbb{N}$, define $$\mathcal{R}_{<n}=\left\{f\in \mathcal{R}:\deg(f)<n\right\}=\left\{\sum_{i=0}^na_ix^i:a_i\in \mathbb{F}_q\right\}.$$
Similarly, let $\mathcal{R}_{=n}=\{f\in \mathcal{R}:\deg(f)=n\}$, let $\mathcal{R}_{\leq n}=\mathcal{R}_{<n}\cup \mathcal{R}_{=n}$, and let $\mathcal{R}_{>n}=\mathcal{R}\setminus\mathcal{R}_{\leq n}$. We say that $f\in \mathcal{R}$ is monic if its leading coefficient is $1$ and denote the set of monic polynomials by $\mathcal{R}_{\text{monic}}$. Let $\mathcal{K}$ be the field of fractions of $\mathcal{R}$, and define the absolute value on $\mathcal{K}$ by $\left|\frac{f}{g}\right|=q^{\deg(f)-\deg(g)}$, where $f,g\in \mathcal{R}$ and $g\neq 0$. Then, the completion of $\mathcal{K}$ with respect to $\vert \cdot\vert$ is the field of formal Laurent series over $\mathbb{F}_q$, defined by
$$\mathcal{K}_{\infty}=\left\{\sum_{n=-h}^{\infty}a_nx^{-n}:a_n\in \mathbb{F}_q\right\}.$$
For $f\in \mathcal{K}_{\infty}$ and $r>0$, define $\mathcal{B}(f,r)=\{g\in \mathcal{K}_{\infty}:\vert f-g\vert\leq r\}$. For $m\geq 1$, define the norm on $\mathcal{K}_{\infty}^m$ by $\Vert \mathbf{v}\Vert=\max_{i=1,\dots,m}\vert v_i\vert$, where $\mathbf{v}=(v_1,\dots,v_m)$. Define analogously $$\mathcal{B}(\mathbf{v},r)=\{\mathbf{u}\in \mathcal{K}_{\infty}^m:\Vert \mathbf{v}-\mathbf{u}\Vert\leq r\}.$$
For every $m\in \mathbb{N}$, the norm $\Vert \cdot \Vert$ satisfies the ultrametric inequality.
$$\Vert \mathbf{f}+\mathbf{g}\Vert\leq \max\{\Vert \mathbf{f}\Vert,\Vert \mathbf{g}\Vert\}.$$
As a consequence, we have the following result about intersections of balls which is a particular case of \cite[Lemma 6.8]{AK}. For a more comprehensive introduction about function fields, see \cite{A,Ros}. 
\begin{lemma}
\label{lem:BallInt}
    Let $m\in \mathbb{N}$, let $r>0$ and let $\boldsymbol{\alpha}_1,\boldsymbol{\alpha}_2\in \mathcal{K}_{\infty}^m$. 
    \begin{enumerate}
        \item If $\Vert \boldsymbol{\alpha}_1-\boldsymbol{\alpha}_2\Vert>r$, then, $\mathcal{B}(\boldsymbol{\alpha}_1,r)\cap \mathcal{B}(\boldsymbol{\alpha}_2,r)=\emptyset$.
        \item If $\Vert \boldsymbol{\alpha}_1-\boldsymbol{\alpha}_2\Vert\leq r$, then, $\mathcal{B}(\boldsymbol{\alpha}_1,r)=\mathcal{B}(\boldsymbol{\alpha}_2,r)$.
    \end{enumerate}
\end{lemma}
Lemma \ref{lem:BallInt} stands at the heart of many of our proofs, and therefore, our methods cannot be trivially generalized to the real setting. Let $\mathcal{O}=\{\alpha\in \mathcal{K}_{\infty}:\vert \alpha\vert\leq 1\}$, and let 
$$\mathfrak{m}=x^{-1}\mathcal{O}=\{\alpha\in \mathcal{K}_{\infty}:\vert \alpha\vert\leq q^{-1}\}.$$
For $\boldsymbol{\alpha}\in \mathcal{K}_{\infty}^m$, we write $\boldsymbol{\alpha}=[\boldsymbol{\alpha}]+\langle \boldsymbol{\alpha}\rangle$, where $[\boldsymbol{\alpha}]\in \mathcal{R}^m$ and $\langle\boldsymbol{\alpha}\rangle\in \mathfrak{m}^m$. In this paper, the Haar measure on $\mathcal{K}_{\infty}$ is defined as the unique translation invariant measure $\nu$, such that $\nu(\mathfrak{m})=1$. For $m\in \mathbb{N}$, let $\nu_m$ denote the product measure $\underbrace{\nu\times \dots\times\nu}_{m\text{ times}}$ on $\mathcal{K}_{\infty}^m$. 
\subsection{Main Results}
Let $\mathcal{S}\subseteq \mathcal{R}_{\text{monic}}$ be an infinite set of polynomials. We say that $\mathcal{S}$ is \emph{degree covering} if $\mathcal{S}\cap \mathcal{R}_{=n}\neq \emptyset$ for every $n\in \mathbb{N}$. For $\mathcal{S}\subseteq \mathcal{R}_{\text{monic}}$ and $n\in \mathbb{N}$, define 
$$m_{\mathcal{S}}(n)=\min\{k\geq n:\mathcal{S}\cap \mathcal{R}_{=k}\neq \emptyset\}.$$
For $m\in \mathbb{N}$, we say that a vector $\mathbf{v}=(v_1,\dots,v_m)\in \mathcal{R}^m$ is \emph{primitive} if $\mathbf{v}\neq a\mathbf{u}$ for $\mathbf{u}\in \mathcal{R}^m$ and $a\in \mathcal{R}\setminus \mathbb{F}_q^*$. Equivalently, $\mathbf{v}$ is primitive if and only if $\gcd(v_1,\dots,v_m)=1$. The set of primitive vectors  in $\mathcal{R}^m$ is denoted by $\widehat{\mathcal{R}}^m$. 
\begin{definition}
For $\boldsymbol{\alpha}\in \mathcal{K}_{\infty}^m$ and $n\in \mathbb{N}$, define
$$\deg_{\min,\mathcal{S}}(\boldsymbol{\alpha},q^{-n})=\min\left\{d:\exists \frac{\mathbf{P}}{Q}\in \mathcal{K}^m:(\mathbf{P},Q)\in \widehat{\mathcal{R}}^{m+1},Q\in \mathcal{S}\cap \mathcal{R}_{=d}, \Vert \mathbf{P}\Vert<\vert Q\vert, \left\Vert\boldsymbol{\alpha}-\frac{\mathbf{P}}{Q}\right\Vert<\frac{1}{q^n}\right\}.$$
\end{definition}
First, observe that for every $\mathbf{a}\in \mathcal{R}^m$ and for every $\boldsymbol{\alpha}\in \mathfrak{m}^m$, we have $\deg_{\min,\mathcal{S}}(\mathbf{a}+\boldsymbol{\alpha},q^{-n})=\deg_{\min,\mathcal{S}}(\boldsymbol{\alpha},q^{-n})$. Hence, it suffices to study the function $\deg_{\min,\mathcal{S}}(\cdot,q^{-n})$ on $\mathfrak{m}^m$. Moreover, by Dirichlet's approximation theorem \cite[Theorem 1.1]{GaGh17} (see also \cite[Theorem 4.1]{A}), for every $\mathcal{S}\subseteq \mathcal{R}_{\text{monic}}$, for every $n\in \mathbb{N}$, and for every $\boldsymbol{\alpha}\in \mathcal{K}_{\infty}^m$, we have $\deg_{\min,\mathcal{S}}(\boldsymbol{\alpha},q^{-n})\leq m_{\mathcal{S}}(n)$. In particular, if $\mathcal{S}$ is degree covering, then for every $\boldsymbol{\alpha}\in \mathcal{K}_{\infty}^m$, we have $\deg_{\min,\mathcal{S}}(\boldsymbol{\alpha},q^{-n})\leq n$. 

In \cite{Ara}, the distribution of $\deg_{\min,\mathcal{R}_{\text{monic}}}(\alpha,q^{-n})$ was computed when $m=1$ by using linear algebra and basic number theory. One can also study the discrete analogue in $\mathcal{K}_{\infty}$. Let $N\in \mathcal{R}$ be a polynomial with $\deg(N)=n$. For $\mathbf{a}\in \mathcal{R}_{<n}^m$, define $d_{N,\mathcal{S}}(\mathbf{a})=\deg_{\min,\mathcal{S}}\left(\frac{\mathbf{a}}{N},\frac{1}{\vert N\vert}\right)$. We prove that $\deg_{\min,\mathcal{S}}(\cdot,\vert N\vert^{-1})$ and $d_{N,\mathcal{S}}(\cdot)$ have the same distribution for every $N\in \mathcal{R}$, for every dimension $m\in \mathbb{N}$, and for every infinite set $\mathcal{S}\subseteq \mathcal{R}_{\text{monic}}$. 
\begin{theorem}
\label{thm:SameDist}
    Let $\mathcal{S}\subseteq \mathcal{R}_{\text{monic}}$ be an infinite set. Then, for every dimension $m\in \mathbb{N}$, for every $n\in \mathbb{N}$, for every $N\in \mathcal{R}_{=n}$, and for every $k=0,\dots,m_{\mathcal{S}}(n)$, we have
    \begin{equation}
        \mathbb{P}\left(\mathbf{a}\in \mathcal{R}_{<n}^m:d_{N,\mathcal{S}}(\mathbf{a})=k\right)=\nu_m\left(\boldsymbol{\alpha}\in \mathfrak{m}^m:\deg_{\min,\mathcal{S}}(\boldsymbol{\alpha},q^{-n})=k\right),
    \end{equation}
    where $\mathbb{P}$ is the uniform probability on $\mathcal{R}_{<n}^m$.
\end{theorem}
As a corollary, the averages of the discrete and continuous minimal denominator functions are equal to one another, which can be viewed as a function field analogue of \cite{BM,Sh}.
\begin{corollary}
    Let $\mathcal{S}\subseteq \mathcal{R}_{\text{monic}}$ and let $m\in \mathbb{N}$. Then, for every $N\in \mathcal{R}$, 
    $$\mathbb{E}[d_{N,\mathcal{S}}(a)]=\frac{1}{\vert N\vert}\sum_{\mathbf{a}\in \mathcal{R}^m}d_{N,\mathcal{S}}(\mathbf{a})=\int_{\mathfrak{m}^n}\deg_{\min,\mathcal{S}}(\boldsymbol{\alpha},\vert N\vert^{-1})d\nu_m(\boldsymbol{\alpha}).$$
\end{corollary}
The following result is an immediate corollary of Theorem \ref{thm:SameDist} and \cite[Theorem 1.1]{Ara}.
\begin{corollary}
    Let $n\in \mathbb{N}$ and let $N\in \mathcal{R}_{=n}$. If $n=1$, we have
    $$\mathbb{P}(a\in \mathcal{R}_{<1}:d_{N,\mathcal{R}_{\text{monic}}}(a)=k)=\begin{cases}
        \frac{1}{q}&k=0\\
        \frac{q-1}{q}&k=1
    \end{cases}$$
    $$=\nu\left(\alpha\in \mathfrak{m}:\deg_{\min,\mathcal{R}_{\text{monic}}}(\alpha,q^{-1})=k\right).$$
    If $n\geq 2$, then, 
    $$\mathbb{P}(a\in \mathcal{R}_{<n}:d_{N,\mathcal{R}_{\text{monic}}}(a)=k)=\begin{cases}
        q^{-n}&k=0\\
        \frac{q-1}{q^{n-2k+1}}&k\leq \left\lceil\frac{n}{2}\right\rceil\\
        0&\text{else}
    \end{cases}$$
    $$=\nu\left(\alpha\in \mathfrak{m}:\deg_{\min,\mathcal{R}_{\text{monic}}}(\alpha,q^{-n})=k\right).$$
\end{corollary}
Moreover, in analogy with \cite{Ara}, one can discuss $\mathcal{S}$ minimal denominators.
\begin{definition}
    Let $m\in \mathbb{N}$, let $n\in \mathbb{N}$, and let $\boldsymbol{\alpha}\in \mathfrak{m}^m$. For an infinite set $\mathcal{S}\subseteq \mathcal{R}_{\text{monic}}$, we say that $Q\in \mathcal{S}$ is an $\mathcal{S}$ minimal denominator for $(\boldsymbol{\alpha},q^{-n})$ if 
    \begin{enumerate}
        \item $\deg(Q)=\deg_{\min,\mathcal{S}}(\boldsymbol{\alpha},q^{-n})$, and
        \item there exists $\mathbf{P}\in \mathcal{R}^m$ such that 
        \begin{enumerate}
            \item $(\mathbf{P},Q)\in \widehat{\mathcal{R}}^{m+1}$, and
            \item $\left|\boldsymbol{\alpha}-\frac{\mathbf{P}}{Q}\right|<q^{-n}$. 
        \end{enumerate}
    \end{enumerate}
\end{definition}
By \cite[Lemma 1.3]{Ara}, for every $\boldsymbol{\alpha}\in \mathfrak{m}^m$ and $n\in \mathbb{N}$, there exists a unique monic $\mathcal{S}$ minimal denominator for $(\boldsymbol{\alpha},q^{-n})$, which we denote by $Q_{\min,\mathcal{S}}(\boldsymbol{\alpha},q^{-n})$. In \cite{Ara}, the distribution of $Q_{\min,\mathcal{R}_{\text{monic}}}(\cdot,q^{-n})$ was computed for every $n\in \mathbb{N}$ when $m=1$. For $N\in \mathcal{R}_{=n}$ and for $\mathbf{a}\in \mathcal{R}_{<n}^m$,  define $Q_{N,\mathcal{S}}(\mathbf{a}):=Q_{\min,\mathcal{S}}\left(\frac{\mathbf{a}}{N},\frac{1}{\vert N\vert}\right)$. 
\begin{theorem}
\label{thm:Q_minEquals}
    Let $\mathcal{S}\subseteq \mathcal{R}_{\text{monic}}$ be an infinite set. Then, for every dimension $m\in \mathbb{N}$, for every $n\in \mathbb{N}$, for every $N\in \mathcal{R}_{=n}$, and for every $Q\in \mathcal{S}$, we have
    \begin{equation}
        \mathbb{P}\left(\mathbf{a}\in \mathcal{R}_{<n}^m:Q_{N,\mathcal{S}}(\mathbf{a})=Q\right)=\nu\left(\boldsymbol{\alpha}\in \mathfrak{m}^m:Q_{\min,\mathcal{S}}(\boldsymbol{\alpha},q^{-n})=Q\right).
    \end{equation}
\end{theorem}
\begin{remark}
    The proofs of Theorem \ref{thm:SameDist} and Theorem \ref{thm:Q_minEquals} rely on a function field analogue of "gaps" for $\mathcal{S}$-Farey fractions with restricted denominators, in a similar fashion, as was done in \cite{Mar,M}.
\end{remark}
In particular, by \cite[Theorem 1.7]{Ara} and Theorem \ref{thm:Q_minEquals}, we obtain the distribution of $Q_{N,\mathcal{R}_{\text{monic}}}$ when $m=1$.
\begin{corollary}
    Let $n\in \mathbb{N}$, and let $N\in \mathcal{R}_{=n}$. Then, for every $Q\in \mathcal{R}_{\leq n}\cap \mathcal{R}_{\text{monic}}$, we have
    $$\mathbb{P}\left(a\in \mathcal{R}_{<n}:Q_{N,\mathcal{R}_{\text{monic}}}(a)=Q\right)$$
    $$=\frac{1}{q^n}\left(\vert Q\vert+\sum_{M\in S(Q)}\vert M\vert\sum_{\ell=1}^{D(M)}(-1)^{\ell}\left(\frac{D\left(\frac{Q}{M}\right)!}{\left(D\left(\frac{Q}{N}\right)-\ell\right)!}+\sum_{R\in S\left(\frac{Q}{M}\right):D\left(\frac{Q}{RM}\right)\geq \ell}\mu(R)\frac{D(R)!}{(D(R)-\ell)!}\right)\right),$$
    where $S(Q)$ is the set of monic divisors of $Q$, $D(Q)=\#S(Q)$, and $\mu(Q)$ is the Mobius function. 
\end{corollary}
Moreover, we use Theorem \ref{thm:Q_minEquals} to compute the distribution of $\deg_{\min,\mathcal{S}}$ and $Q_{\min,\mathcal{S}}$, where $\mathcal{S}$ is a multiplicative semigroup generated by one element, which is often referred to as Hadamard lacunary sequences. These sequences are of the form $\{P^d:d\geq 0\}$ and they relate to the base $P$ expansion of formal Laurent series. 
\begin{theorem}
\label{thm:Lacunary}
    Let $P\in \mathcal{R}\setminus \mathbb{F}_q$, let $m,n\in \mathbb{N}$, and let $\mathcal{S}=\{P^d:d\geq 0\}$. Then, 
    \begin{align}
    \begin{split}
    \label{eqn:P^nDist}
        \mathbb{P}\left(\boldsymbol{\alpha}\in \mathfrak{m}^m:Q_{\min,\mathcal{S}}(\boldsymbol{\alpha},q^{-n})=P^d\right)&=\mathbb{P}\left(\boldsymbol{\alpha}\in \mathfrak{m}^m:\deg_{\min,\mathcal{S}}(\boldsymbol{\alpha},q^{-n})=d\deg(P)\right)\\&=\begin{cases}
            \frac{1}{q^{mn}}&d=0\\
            \frac{\vert P\vert^{md}-\vert P\vert^{m(d-1)}}{q^{mn}}&1\leq d\leq k\\
            \frac{q^{mn}-\vert P\vert^{mk}}{q^{mn}}&d=k+1\\
            0&\text{else}
        \end{cases}.
    \end{split}
    \end{align}
\end{theorem} 
\section{Proof of Theorem \ref{thm:SameDist}}
For $\mathcal{S}\subseteq \mathcal{R}$, $m\in \mathbb{N}$, and $k\in \mathbb{N}$, define the $\mathcal{S}$ Farey fractions of degree at most $k$ by 
$$\mathcal{F}^m_{k,\mathcal{S}}=\left\{\frac{\mathbf{P}}{Q}\in \mathcal{K}:(\mathbf{P},Q)\in \widehat{\mathcal{R}}^{m+1},\Vert \mathbf{P}\Vert<\vert Q\vert\leq q^k,Q\in \mathcal{S}\right\},$$
where $\boldsymbol{0}$ is viewed as the Farey fraction $\frac{\boldsymbol{0}}{1}$. We proof Theorem \ref{thm:SameDist} by explicitly connecting between number of disjoint balls of fixed radius centered at $\mathcal{F}^m_{k,\mathcal{S}}$ and the distribution of $d_{N,\mathcal{S}}$ and of $\deg_{\min,\mathcal{S}}(\cdot,q^{-n})$. For $m,n\in \mathbb{N}$ and $0\leq k\leq n$, let $f_{m,n,\mathcal{S}}(k)$ denote the number of disjoint balls in the union $\mathcal{B}\left(\mathcal{F}^m_{k,\mathcal{S}},q^{-(n+1)}\right):=\bigcup_{\frac{\mathbf{P}}{Q}\in \mathcal{F}^m_{k,\mathcal{S}}}\mathcal{B}\left(\frac{\mathbf{P}}{Q},q^{-(n+1)}\right)$. For every $n\in \mathbb{N}$, we have $\mathcal{B}\left(\mathcal{F}^m_{k,\mathcal{S}},q^{-(n+1)}\right)\subseteq \mathfrak{m}^m$, and moreover, $\mathcal{B}\left(\mathcal{F}_{m_{\mathcal{S}}(n),\mathcal{S}}^m,q^{-(n+1)}\right)=\mathfrak{m}^m$.
\begin{theorem}
    \label{thm:DistDeg}
    Let $\mathcal{S}\subseteq \mathcal{R}_{\text{monic}}$ be an infinite set, let $m,n\in \mathbb{N}$, let $N\in \mathcal{R}_{=n}$, and let $1\leq k\leq n$. Then,
    \begin{equation*}
    \begin{split}
        \nu_m\left(\boldsymbol{\alpha}\in \mathfrak{m}^m:\deg_{\min,\mathcal{S}}(\alpha,q^{-n})=k\right)=\frac{f_{m,n,\mathcal{S}}(k)-f_{m,n,\mathcal{S}}(k-1)}{q^{mn}}=\mathbb{P}\left(\mathbf{a}\in \mathcal{R}^m_{<n}:d_{N,\mathcal{S}}(\mathbf{a})=k\right),
    \end{split}
    \end{equation*}
    and 
    \begin{equation*}
        \nu_m\left(\boldsymbol{\alpha}\in \mathfrak{m}^m:\deg_{\min,\mathcal{S}}(\alpha,q^{-n})=0\right)=\mathbb{P}\left(\mathbf{a}\in \mathcal{R}^m_{<n}:d_{N,\mathcal{S}}(\mathbf{a})=0\right)=\begin{cases}
            q^{-mn}&\text{ if }1\in \mathcal{S}\\
            0&\text{else}
        \end{cases}.
    \end{equation*}
\end{theorem}
Theorem \ref{thm:DistDeg} immediately implies Theorem \ref{thm:SameDist}.
\begin{proof}[Proof of Theorem \ref{thm:DistDeg}]
Let $m,n\in \mathbb{N}$, and let $\boldsymbol{\alpha}\in \mathfrak{m}^m$. For $0\leq k\leq n$, we have $\deg_{\min,\mathcal{S}}(\boldsymbol{\alpha},q^{-n})>k$ if and only if for every $\frac{\mathbf{P}}{Q}\in \mathcal{F}_{k,\mathcal{S}}^m$, we have $\left\Vert\boldsymbol{\alpha}-\frac{\mathbf{P}}{Q}\right\Vert\geq \frac{1}{q^{n}}$. Hence, $\deg_{\min,\mathcal{S}}(\boldsymbol{\alpha},q^{-n})>k$ if and only if $\boldsymbol{\alpha}\in \mathfrak{m}^m\setminus \bigcup_{\frac{\mathbf{P}}{Q}\in \mathcal{F}_{k,\mathcal{S}}^m}\mathcal{B}\left(\frac{\mathbf{P}}{Q},\frac{1}{q^{n+1}}\right)$. Thus, by Lemma \ref{lem:BallInt},
\begin{equation}
    \nu_m\left(\boldsymbol{\alpha}\in \mathfrak{m}^m:\deg_{\min,\mathcal{S}}(\alpha,q^{-n})>k\right)=1-\nu_m\left(\bigcup_{\frac{\mathbf{P}}{Q}\in \mathcal{F}_{k,\mathcal{S}}^m}\mathcal{B}\left(\frac{\mathbf{P}}{Q},\frac{1}{q^{n+1}}\right)\right)=1-\frac{f_{m,n,\mathcal{S}}(k)}{q^{mn}}.
\end{equation}
Thus, for every $k>0$, we have
\begin{equation}
    \nu_m\left(\boldsymbol{\alpha}\in \mathfrak{m}^m:\deg_{\min,\mathcal{S}}(\boldsymbol{\alpha},q^{-n})=k\right)=\frac{f_{m,n,\mathcal{S}}(k)-f_{m,n,\mathcal{S}}(k-1)}{q^{mn}}.
\end{equation}
We now compute the distribution of $d_{N,\mathcal{S}}(\mathbf{a})$. For $0\leq k\leq m_{\mathcal{S}}(n)-1$, we have $d_{N,\mathcal{S}}(\mathbf{a})>k$ if and only if for every $\frac{\mathbf{P}}{Q}\in \mathcal{F}^m_{k,\mathcal{S}}$, we have $\left\Vert\frac{\mathbf{a}}{N}-\frac{\mathbf{P}}{Q}\right\Vert\geq q^{-n}$. Thus, $\left\Vert\mathbf{a}-\frac{N\mathbf{P}}{Q}\right\Vert\geq 1$, so that 
$$d_{N,\mathcal{S}}(\mathbf{a})>k\Leftrightarrow
\mathbf{a}\in \mathcal{R}_{<n}^m\setminus \bigcup_{\frac{\mathbf{P}}{Q}\in \mathcal{F}^m_{k,\mathcal{S}}}\left(\mathcal{B}\left(\frac{N\mathbf{P}}{Q},q^{-1}\right)\cap \mathcal{R}_{<n}^m\right).$$
Since $\left\Vert\frac{N\mathbf{P}}{Q}\right\Vert<q^n$, then, $\#\mathcal{B}\left(\frac{N\mathbf{P}}{Q},q^{-1}\right)\cap \mathcal{R}_{<n}=1$ for every $\frac{\mathbf{P}}{Q}\in \mathcal{F}_{k,\mathcal{S}}^m$. 

Furthermore, by Lemma \ref{lem:BallInt}, $\mathcal{B}\left(\frac{\mathbf{P}}{Q},q^{-(n+1)}\right)\cap \mathcal{B}\left(\frac{\mathbf{A}}{B},q^{-(n+1)}\right)=\emptyset$ if and only if \begin{equation}\label{eqn:BallDiff} \left\Vert\frac{\mathbf{P}}{Q}-\frac{\mathbf{A}}{B}\right\Vert<q^{-n}.\end{equation} Therefore, by multiplying equation \eqref{eqn:BallDiff} by $N$ and again applying Lemma \ref{lem:BallInt}, we have $\mathcal{B}\left(\frac{\mathbf{P}}{Q},q^{-(n+1)}\right)\cap \mathcal{B}\left(\frac{\mathbf{A}}{B},q^{-(n+1)}\right)=\emptyset$ if and only if $\mathcal{B}\left(\frac{N\mathbf{P}}{Q},q^{-1}\right)\cap \mathcal{B}\left(\frac{N\mathbf{A}}{B},q^{-1}\right)=\emptyset$. Thus, the number of disjoint balls in the set $\left\{\mathcal{B}\left(\frac{N\mathbf{P}}{Q},q^{-1}\right):\frac{\mathbf{P}}{Q}\in \mathcal{F}_{k,\mathcal{S}}^m\right\}$ is equal to $f_{m,n,\mathcal{S}}(k)$. Therefore, 
\begin{equation}
    \#\left\{\mathbf{a}\in \mathcal{R}_{<n}^m:d_{N,\mathcal{S}}(\mathbf{a})>k\right\}=\#\mathcal{R}_{<n}^m-f_{m,n,\mathcal{S}}(k)=q^{mn}-f_{m,n,\mathcal{S}}(k).
\end{equation}
As a consequence, for every $k\geq 1$,
\begin{equation}
    \#\{\mathbf{a}\in \mathcal{R}_{<n}^m:d_{N,\mathcal{S}}(\mathbf{a})=k\}=f_{m,n,\mathcal{S}}(k)-f_{m,n,\mathcal{S}}(k-1),
\end{equation}
so that for every $k=1,\dots, m_{\mathcal{S}}(n)$, 
\begin{equation}
\label{eqn:Distwithf_n,S(k)}
\begin{split}
    \mathbb{P}\left(\mathbf{a}\in \mathcal{R}_{<n}^m:d_{N,\mathcal{S}}(\mathbf{a})=k\right)=\frac{f_{m,n,\mathcal{S}}(k)-f_{m,n,\mathcal{S}}(k-1)}{q^{mn}}\\
    =\nu_m\left(\boldsymbol{\alpha}\in \mathfrak{m}^m:\deg_{\min,\mathcal{S}}(\boldsymbol{\alpha},q^{-n})=k\right).
\end{split}
\end{equation}
When $k=0$, the proof is split into two cases. If $1\in \mathcal{S}$, then $\deg_{\min,\mathcal{S}}(\boldsymbol{\alpha},q^{-n})=0$ if and only if $\Vert \boldsymbol{\alpha}\Vert<q^{-n}$. Hence, 
$$\nu_m(\boldsymbol{\alpha}\in \mathfrak{m}^m:\deg_{\min,\mathcal{S}}(\boldsymbol{\alpha},q^{-n})=0)=q^{-nm}.$$
On the other hand, let $N\in \mathcal{R}_{=n}$ and $\mathbf{a}\in\mathcal{R}_{<n}^m$. Note that $d_{N,\mathcal{S}}(\mathbf{a})=0$ if and only if $\mathbf{a}=0$, so that 
$$\mathbb{P}(\mathbf{a}\in \mathcal{R}_{<n}^m:d_{N,\mathcal{S}}(\mathbf{a})=0)=\frac{1}{q^{mn}}=\nu_m\left(\boldsymbol{\alpha}\in \mathfrak{m}^m:\deg_{\min,\mathcal{S}}(\boldsymbol{\alpha},q^{-n})=0\right).$$
If $1\notin \mathcal{S}$, then, $\deg_{\min,\mathcal{S}}(\boldsymbol{\alpha},q^{-n})>0$ for every $n\in \mathbb{N}$ and for every $\boldsymbol{\alpha}\in \mathfrak{m}^m$. In particular $d_{N,\mathcal{S}}(\mathbf{a})=\deg_{\min,\mathcal{S}}\left(\frac{\mathbf{a}}{N},q^{-n}\right)>0$ for every $\mathbf{a}\in \mathcal{R}^{m}_{<n}$. Hence, 
$$\mathbb{P}(\mathbf{a}\in \mathcal{R}_{<n}^m:d_{N,\mathcal{S}}(\mathbf{a})=0)=0=\nu_m(\boldsymbol{\alpha}\in\mathfrak{m}^m:\deg_{\min,\mathcal{S}}(\boldsymbol{\alpha},q^{-n})=0).$$
\end{proof}
Theorem \ref{thm:DistDeg} gives rise to the following corollary. 
\begin{corollary}
\label{cor:d_minFarey}
    Let $\mathcal{S}\subseteq \mathcal{R}_{\text{monic}}$, let $m,n\in \mathbb{N}$, and let $N\in \mathcal{R}_{=n}$. Then, for every $1\leq k\leq \frac{n}{2}$, we have 
    \begin{equation}\label{eqn::DistSmallk}\nu\left(\boldsymbol{\alpha}\in \mathfrak{m}^m:\deg_{\min,\mathcal{S}}(\boldsymbol{\alpha},q^{-n})=k\right)=\frac{\#\mathcal{F}^m_{k,\mathcal{S}}-\#\mathcal{F}^m_{k-1,\mathcal{S}}}{q^{mn}}=\mathbb{P}\left(\mathbf{a}\in \mathcal{R}_{<n}^m:d_{N,\mathcal{S}}(\mathbf{a})=k\right).\end{equation}
\end{corollary}
\begin{proof}
    By Theorem \ref{thm:DistDeg}, it suffices to compute $f_{m,n,\mathcal{S}}(k)-f_{m,n,\mathcal{S}}(k-1)$. If $\frac{\mathbf{A}}{B},\frac{\mathbf{P}}{Q}\in \mathcal{F}^m_{k,\mathcal{S}}$ satisfy $\left\Vert\frac{\mathbf{A}}{B}-\frac{\mathbf{P}}{Q}\right\Vert<\frac{1}{q^n}$, then, 
    $$\frac{1}{q^n}>\left\Vert\frac{Q\mathbf{A}-B\mathbf{P}}{BQ}\right\Vert\geq \frac{1}{\vert BQ\vert}>\frac{1}{q^{2k}}.$$
    Thus, $2k>n$. Therefore, for every $k\leq \frac{n}{2}$, the set $\mathcal{B}(\mathcal{F}^m_{k,\mathcal{S}},q^{-(n+1)})$ is composed of disjoint balls of radius $q^{-(n+1)}$. Therefore, $f_{m,n,\mathcal{S}}(k)=\#\mathcal{F}^m_{k-1,\mathcal{S}}$, so that \eqref{eqn::DistSmallk} follows.  
\end{proof}
Due to Theorem \ref{thm:SameDist} and Corollary \ref{cor:d_minFarey}, it is natural to ask the following question.
\begin{question}
    Let $\mathcal{S}\subseteq \mathcal{R}_{\text{monic}}$ be an infinite set. For $k,m,n\in \mathbb{N}$, what is $f_{m,n,\mathcal{S}}(k)$?
\end{question}
\section{Proof of Theorem \ref{thm:Q_minEquals}}
To prove Theorem \ref{thm:Q_minEquals}, we first observe the following property of minimal denominators. 
\begin{lemma}
\label{lem:SepDenoms}
    Let $\mathcal{S}\subseteq \mathcal{R}_{\text{monic}}$ be an infinite set, let $m,n\in \mathbb{N}$, let $k\leq m_{\mathcal{S}}(n)$, let $Q\in \mathcal{R}_{=k}\cap \mathcal{S}$, and let $\boldsymbol{\alpha}\in \mathfrak{m}^m$. Then, $Q_{\min,\mathcal{S}}(\boldsymbol{\alpha},q^{-n})=Q$ implies that there exists some $\mathbf{P}\in \mathcal{R}_{<k}^m$ such that $(\mathbf{P},Q)\in \widehat{\mathcal{R}}^{m+1}$, and for every $\frac{\mathbf{A}}{B}\in \mathcal{F}^m_{k-1,\mathcal{S}}$, one has $\left\Vert\frac{\mathbf{A}}{B}-\frac{\mathbf{P}}{Q}\right\Vert\geq q^{-n}$.
\end{lemma}
\begin{proof}
    Note that $Q_{\min,\mathcal{S}}(\boldsymbol{\alpha},q^{-n})=Q$ if and only if 
    $$\boldsymbol{\alpha}\in \bigcup_{\begin{matrix}\mathbf{P}\in \mathcal{R}_{<k}^m\\ 
    (\mathbf{P},Q)\in \widehat{\mathcal{R}}^{m+1}\end{matrix}}\mathcal{B}\left(\frac{\mathbf{P}}{Q},q^{-(n+1)}\right)\setminus \bigcup_{\frac{\mathbf{A}}{B}\in \mathcal{F}^m_{k-1,\mathcal{S}}}\mathcal{B}\left(\frac{\mathbf{A}}{B},q^{-(n+1)}\right).$$
    By Lemma \ref{lem:BallInt}, the balls $\mathcal{B}\left(\frac{\mathbf{P}}{Q},q^{-(n+1)}\right)$ and $\mathcal{B}\left(\frac{\mathbf{A}}{B},q^{-(n+1)}\right)$ are either disjoint or equal to one another. Hence, if $Q_{\min,\mathcal{S}}(\boldsymbol{\alpha},q^{-n})=Q$, then, there exists $\mathbf{P}\in \mathcal{R}_{<k}^m$ with $(\mathbf{P},Q)\in \widehat{\mathcal{R}}^{m+1}$, such that for every $\frac{\mathbf{A}}{B}\in \mathcal{F}_{k-1,\mathcal{S}}^m$, we have $\mathcal{B}\left(\frac{\mathbf{P}}{Q},q^{-(n+1)}\right)\cap \mathcal{B}\left(\frac{\mathbf{A}}{B},q^{-(n+1)}\right)=\emptyset$. Therefore, by Lemma \ref{lem:BallInt}, $\left\Vert\frac{\mathbf{P}}{Q}-\frac{\mathbf{A}}{B}\right\Vert\geq q^{-n}$ for every $\frac{\mathbf{A}}{B}\in \mathcal{F}_{k-1,\mathcal{S}}^m$. 
\end{proof}
This motivates the following definition.
\begin{definition}
    A Farey fraction $\frac{\mathbf{P}}{Q}\in \mathcal{F}_{k,\mathcal{S}}^m$ is called $(\mathcal{S},n)$ separated if for every $\frac{\mathbf{A}}{B}\in \mathcal{F}_{\deg(Q)-1,\mathcal{S}}^m$, we have $\left\Vert\frac{\mathbf{P}}{Q}-\frac{\mathbf{A}}{B}\right\Vert\geq q^{-n}$. For $k\leq m_{\mathcal{S}}(n)$, let $f_{m,n,\mathcal{S},Q}(k)$ denote the number of distinct balls in the union 
    \begin{equation}\label{eqn:B(P,Q)Sep}\bigcup_{\begin{matrix}\mathbf{P}\in \mathcal{R}_{<k}^m\\ \frac{\mathbf{P}}{Q}\in \mathcal{F}_{k,\mathcal{S}}^m\text{ is }(\mathcal{S},n)\text{separated}\end{matrix}}\mathcal{B}\left(\frac{\mathbf{P}}{Q},q^{-(n+1)}\right).\end{equation}
\end{definition}
To prove Theorem \ref{thm:Q_minEquals}, we prove the following theorem which explicitly connects between $f_{m,n,\mathcal{S},Q}(k)$ and the distributions of $Q_{\min,\mathcal{S}}$ and $Q_{N,\mathcal{S}}$.
\begin{theorem}
\label{thm:Q_minDist}
Let $\mathcal{S}\subseteq \mathcal{R}_{\text{monic}}$ be an infinite set, let $m,n\in \mathbb{N}$, let $k\leq m_{\mathcal{S}}(n)$, let $N\in \mathcal{R}_{=n}$, and let $Q\in \mathcal{R}_{=k}\cap \mathcal{S}$. Then,
\begin{equation}
    \nu_m\left(\boldsymbol{\alpha}\in \mathfrak{m}^m:Q_{\min,\mathcal{S}}(\boldsymbol{\alpha},q^{-n})=Q\right)=\frac{f_{m,n,\mathcal{S}}(k)}{q^{{mn}}}=\mathbb{P}\left(\mathbf{a}\in \mathcal{R}_{<n}^m:Q_{N,\mathcal{S}}(\mathbf{a})=Q\right)
\end{equation}
\end{theorem}
Theorem \ref{thm:Q_minEquals} arises as an immediate corollary of Theorem \ref{thm:Q_minDist}.
\begin{proof}[Proof of Theorem \ref{thm:Q_minDist}]
    Let $k\leq m_{\mathcal{S}}(n)$, and let $Q\in \mathcal{R}_{=k}\cap \mathcal{S}$ be a monic polynomial.
    Then, by Lemma \ref{lem:BallInt} and Lemma \ref{lem:SepDenoms},
    \begin{align*}
        \nu_m\left(\boldsymbol{\alpha}\in \mathfrak{m}^m:Q_{\min,\mathcal{S}}(\boldsymbol{\alpha},q^{-n})=Q\right)=\nu_m\left(\bigcup_{\begin{matrix}
            \mathbf{P}\in \mathcal{R}_{<k}^m\\
            \frac{\mathbf{P}}{Q}\text{ is }(\mathcal{S},n)\text{ separated}
        \end{matrix}}\mathcal{B}\left(\frac{\mathbf{P}}{Q},q^{-(n+1)}\right)\right)=\frac{f_{m,n,\mathcal{S},Q}(k)}{q^{mn}}.
    \end{align*}
    On the other hand, if $N\in \mathcal{R}_{=n}$ and $\mathbf{a}\in \mathcal{R}_{<n}^m$, then, by Lemma \ref{lem:SepDenoms}, $Q_{N,\mathcal{S}}(\mathbf{a})=Q$ if and only if there exists $\mathbf{P}\in \mathcal{R}_{<k}^m$ such that $(\mathbf{P},Q)\in \widehat{\mathcal{R}}^{m+1}$, $\frac{\mathbf{P}}{Q}$ is $(\mathcal{S},n)$ separated, and $\left\Vert\frac{\mathbf{a}}{N}-\frac{\mathbf{P}}{Q}\right\Vert<q^{-n}$. Thus, 
    \begin{equation}
    \label{eqn:aUnion}
        \mathbf{a}\in \bigcup_{\begin{matrix}\mathbf{P}\in \mathcal{R}_{<k}^m\\
        \frac{\mathbf{P}}{Q}\text{ is }(\mathcal{S},n) \text{ separated}\end{matrix}}\mathcal{B}\left(\frac{N\mathbf{P}}{Q},q^{-1}\right)\cap \mathcal{R}^m_{<n}.
    \end{equation}
    Note that $\mathcal{B}\left(\frac{N\mathbf{P}}{Q},q^{-1}\right)\cap \mathcal{R}_{<n}^m$ contains exactly one point. Hence, the number of disjoint balls in the union on the right hand side of \eqref{eqn:aUnion} is equal to the number of disjoint balls in the union \eqref{eqn:B(P,Q)Sep}, which is equal to $f_{m,n,\mathcal{S},Q}(k)$ by definition. Thus, by Lemma \ref{lem:BallInt}, 
    $$\mathbb{P}\left(\mathbf{a}\in \mathcal{R}_{<n}^m:Q_{N,\mathcal{S}}(\mathbf{a})=Q\right)=\frac{f_{m,n,\mathcal{S},Q}(k)}{q^{mn}}=\nu_m\left(\boldsymbol{\alpha}\in \mathfrak{m}^m:Q_{\min,\mathcal{S}}(\boldsymbol{\alpha},q^{-n})=Q\right).$$
\end{proof}
\begin{remark}
Due to the proof of Theorem \ref{thm:Q_minDist}, it is natural to ask what is the number of $(\mathcal{S},n)$ separated Farey fractions in $\mathcal{F}_{k,\mathcal{S}}^m$, which can be viewed as a function field analogue of the results of \cite{Mar}.
\end{remark}
\section{Minimal Denominators with Respect to Multiplicative Lacunary Semigroups}
\begin{proof}[Proof of Theorem \ref{thm:Lacunary}]
Let $\mathcal{S}=\{P^d:d\geq 0\}$ be a semigroup generated by one element. Note that 
$$\mathcal{S}\cap \mathcal{R}_{=d}=\begin{cases}
    \{P^r\} & r\deg(P)=d\\
    \emptyset & \text{else}
\end{cases},$$ and therefore, computing the distribution of $Q_{\min,\mathcal{S}}(\cdot,q^{-n})$ analogous to computing the distribution of $\deg_{\min,\mathcal{S}}(\cdot ,q^{-n})$. For every $m\in \mathbb{N}$ and for every $\boldsymbol{\alpha}\in \mathfrak{m}^m$, we have $Q_{\min,\mathcal{S}}(\boldsymbol{\alpha},q^{-n})\in \{1,P,\dots,P^k,P^{k+1}\}$, where $k\deg(P)\leq n<(k+1)\deg(P)$. By Theorem \ref{thm:Q_minEquals}, to compute the distribution of $Q_{\min,\mathcal{S}}(\cdot,q^{-n})$ it suffices to compute the distribution of $Q_{x^rP^k,\mathcal{S}}$, where $n=r+k\deg(P)$ for some $r=0,\dots, \deg(P)-1$. Let $\mathbf{A}\in \mathcal{R}_{<n}^m$ and let $d\leq k$. Then, $Q_{x^rP^k,\mathcal{S}}(\mathbf{A})=P^d$ implies that there exists $\mathbf{B}\in \mathcal{R}^m_{<k\deg(P)}$ such that 
\begin{equation}
    \left\Vert\frac{\mathbf{A}}{x^rP^k}-\frac{\mathbf{B}}{P^d}\right\Vert<\frac{1}{q^n}.
\end{equation}
Therefore, $Q_{x^rP^k,\mathcal{S}}(\mathbf{A})=P^d$ if and only if $d$ is the smallest non-negative integer such that
\begin{equation}
    \left\Vert\mathbf{A}-\mathbf{B}x^rP^{k-d}\right\Vert<1. 
\end{equation}
As a consequence $\mathbf{A}\equiv 0\mod x^rP^{k-d}$, so that $\mathbf{A}\in \mathcal{R}^m_{<n}\cap x^rP^{k-d}\mathcal{R}^m=x^rP^{k-d}\mathcal{R}^m_{<d\deg(P)}$. Hence,
\begin{equation}
\begin{split}
    \mathbb{P}\left(\mathbf{A}\in \mathcal{R}^m_{<n}:d_{x^rP^{k-d},\mathcal{S}}(\mathbf{A})=d\deg(P)\right)=\mathbb{P}\left(\mathbf{A}\in \mathcal{R}_{<n}^m:Q_{x^rP^{k-d},\mathcal{S}}(\mathbf{A})=P^d\right)\\
    =\frac{\#x^rP^{k-d}\mathcal{R}_{<d\deg(P)}^m-\#x^rP^{k-d+1}\mathcal{R}^{m}_{<(d-1)\deg(P)}}{q^{nm}}=\frac{\vert P\vert^{md}-\vert P\vert^{m(d-1)}}{q^{nm}}.
\end{split}
\end{equation}
Moreover, $\mathbb{P}\left(\mathbf{A}\in \mathcal{R}^m_{<n}:d_{x^rP^{k-d},\mathcal{S}}(\mathbf{A})=0\right)=\mathbb{P}\left(\mathbf{A}\in \mathcal{R}_{<n}^m:Q_{x^rP^{k-d},\mathcal{S}}=1\right)=\frac{1}{q^{nm}}$ and 
\begin{equation}
\begin{split}
    \mathbb{P}\left(\mathbf{A}\in \mathcal{R}^m_{<n}:d_{x^rP^{k-d},\mathcal{S}}(\mathbf{A})=(k+1)\deg(P)\right)=\mathbb{P}\left(\mathbf{A}\in \mathcal{R}^m_{<n}:Q_{x^rP^{k-d},\mathcal{S}}(\mathbf{A})=P^{k+1}\right)\\
    =\mathbb{P}\left(\mathbf{A}\in \mathcal{R}_{<n}^m:\mathbf{A}\neq 0\mod x^r\right)=\frac{\#\mathcal{R}_{<n}^m-\#x^r\mathcal{R}^{m}_{<n-r}}{q^{nm}}=\frac{q^{nm}-q^{m(n-r)}}{q^{nm}}.
\end{split}
\end{equation}
\end{proof}
\section{Acknowledgments}
I would like to thank Igor Shparlinski for suggesting this question to me, discussions regarding minimal denominators, and providing comments on an earlier version of this paper. I would also like to thank the anonymous referee, whose comments improved the quality of this paper. \footnote{This manuscript has no data and there is no conflict of interest.}
\bibliography{Ref}
\bibliographystyle{amsalpha}
\end{document}